\documentclass[a4paper,14pt]{article}

\usepackage{latexsym}
\usepackage{amsmath}
\usepackage{amsfonts}
\usepackage{amssymb}
\usepackage{cite}
\usepackage{mathrsfs}
\newtheorem{lem}{Lemma}[section]
\newtheorem{pro}[lem]{Proposition}
\newtheorem{thm}[lem]{Theorem}
\newtheorem{rem}[lem]{Remark}
\newtheorem{cor}[lem]{Corollary}

\renewcommand{\a }{\alpha }
\renewcommand{\b }{\beta }

\renewcommand{\d}{\delta }
\newcommand{\D }{\Delta }

\newcommand{\e }{\varepsilon }

\newcommand{\g }{\gamma}

\renewcommand{\l }{\lambda }

\newcommand{\n }{\nabla }
\newcommand{\vp }{\varphi }

\newcommand{\s }{\sigma }

\renewcommand{\phi}{\varphi}

\renewcommand{\O }{\Omega }
\newcommand{\ov}{\overline}

\newcommand{\be}{\begin{equation}}
\newcommand{\ee}{\end{equation}}
\newenvironment{pf}{\noindent{\sc Proof}.\enspace}{\rule{2mm}{2mm}\medskip}

\newcommand{\R}{\mathbb{R}}
\newcommand{\N}{\mathbb{N}}

\newcommand{\de}{\partial}

\newcommand{\ti}{\tilde}

\newcommand{\M}{\mathcal{M}}

\newcommand{\ra}{{\rangle}}
\newcommand{\la}{{\langle}}

\newcommand{\tn}{\tilde \n}
\DeclareMathOperator{\dist}{dist}
\DeclareMathOperator{\cA}{\mathcal{A}}

\title{Sharp local estimates for the Szeg\"o-Weinberger profile in Riemannian manifolds}
\author{Mouhamed Moustapha Fall   and   Tobias Weth}

\begin{document}
\date{}
\maketitle
 \let\thefootnote\relax\footnotetext{mouhamed.m.fall@aims-senegal.org (M. M. Fall), weth@math.uni-frankfurt (T. Weth).}
\let\thefootnote\relax\footnotetext{ African Institute for Mathematical Sciences of Senegal, 
KM 2, Route de Joal, B.P. 14 18. Mbour, S\'en\'egal.}
\let\thefootnote\relax\footnotetext{ Goethe-Universit\"{a}t Frankfurt, Institut f\"{u}r Mathematik.
Robert-Mayer-Str. 10 D-60054 Frankfurt, Germany.}
\bigskip

\begin{abstract}
We study the local Szeg\"o-Weinberger profile in a
geodesic ball $B_g(y_0,r_0)$ centered at a point $y_0$
in a Riemannian manifold $(\M,g)$. This profile is obtained by maximizing
the first nontrivial Neumann eigenvalue $\mu_2$ of the Laplace-Beltrami
Operator $\Delta_g$ on $\M$ among subdomains of $B_g(y_0,r_0)$ with
fixed volume. We derive a sharp asymptotic bounds of this profile
in terms of the scalar curvature of $\M$ at $y_0$. As a corollary, we deduce
a local comparison principle depending only on the scalar curvature.
Our study is related to previous results on the profile corresponding
to the minimization of the first Dirichlet eigenvalue of $\Delta_g$,
but additional difficulties arise due to the fact that $\mu_2$
is degenerate in the unit ball in $\R^N$ and geodesic balls do not
yield the optimal lower bound in the asymptotics we obtain.
\end{abstract}

\section{Introduction}
Let  $(\M,g)$ be a complete Riemannian manifold of dimension
$N$, $N\geq2$. For a bounded regular domain $\O\subset \M$ we consider the  Neumann eigenvalue problem
\begin{equation}
  \label{eq:16}
\D_g f+\mu\,f=0\quad \text{in $\O$},\qquad  \la \n f,\eta\ra_g=0\quad
\text{on $\de\O$},
\end{equation}
where $\D_g f=div_g(\n f)$ is the Laplace-Beltrami operator on $\M$ and $\eta$ is the outer unit normal to $\de\O$.
The set of eigenvalues, counted with multiplicities, in the above
eigenvalue problem is given as an increasing sequence
$$
0=\mu_1(\O,g)<\mu_2(\O,g)\leq\dots+\infty.
$$
By results of Szeg\"{o} \cite{Sz} and Weinberger \cite{Wein},
balls maximize $\mu_2$ among domains having fixed volume in
$\M=\R^{N}$. More precisely, in \cite{Sz} this was proved for the
planar case $N=2$, whereas in \cite{Wein} the case $N \ge 3$ was
considered.
As remarked in \cite{chavel} and \cite{Ashb-Bengu}, this
result extends to the case of the $N$-dimensional hyperbolic space. Moreover, the
same conclusion holds for  domains contained in a hemisphere
\cite{Ashb-Bengu} and -- under further restrictions on the domain -- also 
in rank-1 symmetric spaces \cite{AS}.\\
The aim of the present paper is to study the geometric variational
problem of maximizing $\mu_2(\O,g)$ among domains with fixed volume
{\em locally} in a general complete Riemannian manifold
$(\M,g)$. In order to state our results, we need to introduce some
notations. For a subset $\Omega \subset \M$, we let $|\Omega|_g$ denote
the volume of $\Omega$ with respect to the metric $g$. For $0<v< |\M|_g$, we define the \textit{Szeg\"o-Weinberger profile} of $\M$ as
$$
SW_{\M}(v,g):=\sup_{\O\subset\M,\,|\O|_g=v}\,\mu_2(\O,g).
$$
Here and in the following, we assume without further mention that only
regular bounded domains $\O\subset\M$ are considered. For open subsets $\cA
\subset \M$ and $0<v< |\cA|_g$, we also define 
$$
SW_{\cA}(v,g):=\sup_{\O\subset\cA,\,|\O|_g=v}\,\mu_2(\O,g),
$$
assuming again without further mention that only regular bounded domains $\O
\subset \cA$ are considered. By Weinberger's result in \cite{Wein}, we then have
$$
SW_{\R^{N}}(v)=\left(\frac{|B|}{v}\right)^{\frac{2}{N}}\,\mu_2(B).
$$
where $B$ denotes the unit ball in $\R^{N}$. The
eigenvalue $\mu_2(B)$ has multiplicity
$N$ with corresponding eigenfunctions  $x\mapsto \vp(|x|)\frac{x_i}{|x|}$, $i=1,\dots,N$,
 where $\vp$ can be expressed in terms of a rescaled Bessel function
 of the first kind and satisfies $\vp(0)=\vp'(1)=0$. For matters of
 convenience, we normalize $\vp$ such that 
\be\label{eq:normvp}
\int_0^1\vp^2(t) t^{N-1}dt=\frac{1}{|B|},
\ee
see Section \ref{s:pn} below. We are interested in the local effect
of curvature terms on the {Szeg\"o-Weinberger profile}. For this we study the profile in a small geodesic ball
$B_g(y_0,r)$ of $\M$ centered at a point $y_0 \in \M$ with radius $r$.
In our main result, we obtain the following optimal two-sided local bound.

\begin{thm}\label{pro:Main-rslt}
Let $\M$ be a complete $N$-dimensional Riemannian manifold with $N\geq 2$, and let
$S$ denote the scalar curvature function on $\M$. Moreover, let $y_0
\in \M$, and let  
\begin{align}
  \label{eq:21}
\g_N&=\frac{ 2\mu_2(B)+(N+2)(N-2)-(N+2)|B|\vp^2(1)   }{
6N(N+2)\mu_2(B)}\nonumber\\
&=\frac{1}{3N(N+2)} + \frac{N-2}{6N\mu_2(B)}
-\frac{1}{3N(\mu_2(B)-N+1)}
\end{align}
Then we have:
\begin{itemize}
\item[(i)] As $v \to 0$, 
\begin{equation}
  \label{eq:24}
\frac{SW_{\M}(v)}{SW_{\R^{N}}(v)} \;\ge \;1- \g_N S(y_0) \Bigl( \frac{v}{|B|}
\Bigr)^{\frac{2}{N}} +o(v^{\frac{2}{N}}).
\end{equation}
\item[(ii)] For every $y_0 \in \M$ and every $\e>0$ , there exists $r_\e>0$ such that
\begin{equation}
  \label{eq:5}
\!\!\!\!\!\!\!\!\!\!\!\!\!\!1-\left(\g_N\,S(y_0) +\e \right)
 \Bigl( \frac{v}{|B|} \Bigr)^{\frac{2}{N}} \;\le\;
 \frac{SW_{\!B_g(y_0,r_\e)}(v)}{SW_{\R^{N}}(v)} \;\leq\; 1-\left(\g_N\,S(y_0) -\e \right)
 \Bigl(\frac{v}{|B|} \Bigr)^{\frac{2}{N}}
\end{equation}
for $v\in \bigl(0\,,\,| B_g(y_0,r_\e)|_g\bigr)$.
\item[(iii)] $\gamma_N<0$ for all $N \ge 2$, and $\gamma_N \to 0$ as
  $N \to \infty$.
\end{itemize}
\end{thm}

Some remarks are in order. The right hand side of (\ref{eq:24}) can be replaced by
$$
1- \gamma_N \Bigl[\sup \limits_{y_0 \in \M} S(y_0)\Bigr] \Bigl( \frac{v}{|B|}
\Bigr)^{\frac{2}{N}} +o(v^{\frac{2}{N}})
$$ 
if the supremum of the
scalar curvature is attained
on $\M$, e.g. if $\M$ is compact. The equality
(\ref{eq:21}) is derived from an integral identity for Bessel functions
 which gives $|B| \vp^2(1)= \frac{2\mu_2(B)}{\mu_2(B)-N+1}$, see Lemma
 \ref{sec:prel-notat} below.  The coefficient $\gamma_N$ is uniquely determined by the
two-sided estimate (\ref{eq:5}) and therefore
sharp. To determine the sign of $\gamma_N$ for large $N$, fine
estimates on $\mu_2(B)$ are needed. In Lemma~\ref{sec:prel-notat}
below, we will prove that 
\begin{equation}
  \label{eq:22}
N+1 \le \mu_2(B) < 
\left \{
\begin{aligned}
&N+2,&& \qquad \qquad \text{for $N=2,3,4$,}\\
&N+1 + \frac{2}{N-2}&&\qquad \qquad \text{for $N \ge 5$.}  
\end{aligned}
\right.    
\end{equation}
From this we then deduce part (iii) of
Theorem~\ref{pro:Main-rslt}. The bounds in (\ref{eq:22}) might not be
new, but we could not find any suitable bound in the literature.  
It seems natural to deduce bounds on $\mu_2(B)$ from the fact that $\sqrt{\mu_2(B)}$ is the first positive zero of the derivative of the function $t \mapsto
  t^{(2-N)/2} J_{N/2}(t)$, where $J_{N/2}$ is the Bessel function of
  the first kind of order $N/2$. However, to obtain the upper bound in
  (\ref{eq:22}), we use the variational characterization of $\mu_2(B)$
  instead. See
  Section~\ref{s:pn} below for details.

As a consequence of Theorem~\ref{pro:Main-rslt}, we readily deduce
the following local isochoric comparison principle related to the
Szeg\"o-Weinberger profile. 

\begin{cor}\label{th:main-th-2d}
Let $(\M_1,g_1)$, $(\M_2,g_2)$ be two $N$-dimensional complete Riemannian manifolds,
$N \ge 2$ 
with scalar curvature functions $S_1$, $S_2$ respectively. Let
$y_1\in\M_1$ and $y_2 \in \M_2$ such that $S_1(y_1)<S_2(y_2)$. Then there exists $r>0$ such that 
  \begin{equation}
    \label{eq:26}
SW_{B_{g_1}(y_1,r)}(v) < SW_{B_{g_2}(y_2,r)}(v) 
  \end{equation}
for any $v \in (0,\min\{|B_{g_1}(y_1,r)|_{g_1},|B_{g_2}(y_2,r)|_{g_2}\}).$  
\end{cor}

We emphazise that in the special case where $(\M_2,g_2)$ is a space
form of constant curvature, the right hand side in
(\ref{eq:26}) may be replaced with $\mu_2(E,g_2)$,
where $E$ is any geodesic ball 
  of volume $v$ in $\M_2$. This follows from the local expansion
  of $\mu_2$ in small geodesic balls in these manifolds, see
  Remark~\ref{sec:expansion-mu_2-small}(ii) below.

Corollary \ref{th:main-th-2d} should be seen in comparison with the results in \cite{Druet-Iso,Fall-eigen,Druet-FK} 
concerning the isoperimetric profile $I_\M$ and the Faber-Krahn profile $FK_\M$ of $\M$. More precisely, set
$$
I_{\M}(v,g):=\inf_{\O\subset \M,\,|\O|_g=v }|\de \O|_g
$$
and
$$
FK_{\M}(v,g):=\inf_{\O\subset\M,\,|\O|_g=v}\,\l_1(\O,g),
$$
with $\l_1(\O,g)$ being the first  Dirichlet eigenvalue of $-\Delta_g$
in $\O$. Let $y\in\M$ and $k \in \R$ be such that
${S}(y)<({N-1})N\,k$, where $S(y)$ denotes the scalar curvature of
$\M$ at $y$. Furthermore, let $(\mathbb{M}^{N},g_k)$ denote the
space form of constant sectional curvature $k$. Then there exists
$r_y>0$ such that for any $v\in \left(0\,,\,\left| B_g(y,r_y)\right|_g
\right)$ and any geodesic ball $E$ of volume $v$ in $(\mathbb{M}^{N},g_k)$, we have
\be\label{eq:Ip}
I_{B_g(y,r_y)}(v,g)>|\de E|_{g_k},
\ee
\be\label{eq:FKp}
FK_{B_g(y,r_y)}(v,g)>\l_1( E,{g_k}).
\ee
Inequality \eqref{eq:Ip} was established by Druet\cite{Druet-Iso}, and
\eqref{eq:FKp} was derived independently by Druet\cite{Druet-FK} and
the first author\cite{Fall-eigen}. The first ingredient in the proof of \eqref{eq:FKp} is the following expansion of
$\l_1(B_g(y,r),g)$ when $r \to 0$:
\begin{equation}
  \label{eq:15}
\lambda_1(B_g(y,r),g)= \frac{\lambda_1(B)}{r^2}-\frac{S(y_0)}{6} +O(r)
\end{equation}
This expansion had already been obtained by
Chavel in \cite[Chapter 8]{chavel}. In the proof of Theorem~\ref{pro:Main-rslt}, we
need to derive a corresponding expansion for $\mu_2(B_g(y,r),g)$. This
is more difficult since $\mu_2(B)$ is degenerate with multiplicity $N$
and the corresponding eigenfunctions are nonradial. As a consequence,
an anisotropic curvature term appears in the corresponding expansion. More
precisely, we have 
\begin{equation}
  \label{eq:28}
\mu_2(B_g(y_0,r),g) =\frac{\mu_2(B)}{r^2}+\a_N^-{S}(y_0) +2
\a_N^+R_{min}(y_0)+o(1) \qquad \text{as $r \to 0$}
\end{equation}
with suitable constants $\alpha_N^\pm$ and $R_{min}(y_0)=
\inf\{Ric_{y_0}(A,A) \,:\, A\in T_{y_0}\M,\:|A|=1\}$, see Proposition~\ref{lem:expa_un2B} below. In order to obtain an 
expansion depending only on the scalar curvature, we need to consider
suitable geodesic ellipsoids with small eccentricity. This is a crucial
step in the proof of Theorem~\ref{pro:Main-rslt}, since -- in contrast to
the Faber-Krahn profile -- geodesic balls
do not give rise to optimal two-sided bounds. As a further tool, we need a
quantitative version of the {Szeg\"o-Weinberger} inequality, which has
been obtained very recently in the euclidean case by Brasco and
Pratelli \cite{BP}. In the proof of Theorem~\ref{pro:Main-rslt} we combine these tools with variants of ideas in
\cite{Fall-eigen} and \cite{Wein,Ashb-Bengu,AS} to control error
terms and to construct suitable test functions for the variational
characterization of $\mu_2$, see Section \ref{s:Ubm2} below.
\\
We like to mention that \cite{Fall-eigen} also contains a
statement about the local expansion of a profile related to {\em minimizing}
$\mu_2$ among domains of fixed volume relative to an open set,
see \cite[Theorem 1.3]{Fall-eigen}. However, the proof of this
statement is not correct since it relies on a comparison with a
relative isoperimetric profile which does not correspond to the Neumann
boundary conditions in (\ref{eq:16}) but rather to mixed boundary
conditions.\\
Theorem~\ref{pro:Main-rslt} gives a first hint that critical domains for $\mu_2$ which are nearly balls, if they exists,
 might be located near critical points of the scalar curvature of $\M$
(at least in the twodimensional case). Here, roughly speaking, by a critical
domain we mean a domain where $\mu_2$ is critical with respect to
volume preserving perturbations. Pacard and Sicbaldi \cite{Pac-Sic}
showed that close to nondegenerate critical points of the scalar
curvature there exist small critical domains for the first Dirichlet
eigenvalue of $\Delta_g$. The corresponding
problem for the Neumann eigenvalue $\mu_2$ seems much more
difficult. We note that Zanger \cite{Zanger} derived a Hadamard
type formula (in the spirit of \cite[p. 522]{hadamard}) for a Neumann
eigenvalue which depends smoothly on domain
variations. However, due to possible degeneracy, $\mu_2$ might not
depend smoothly on domain
variations, and therefore it is
not clear how critical domains should be defined. On the other hand, in
\cite{Soufi-Ilias} a notion of critical domains for higher Dirichlet
eigenvalues, which may also be degenerate, is derived via analytic
perturbation theory. It therefore seems natural -- but far from obvious
-- to develop and analyze a similar notion for $\mu_2$. We wish to
address this problem in future work. 

The paper is organized as follows. In Section~\ref{s:pn} we collect
properties of $\mu_2(B)$ and the function $\vp$ appearing in the
definition of the corresponding eigenfunctions. In particular, we
prove the bounds (\ref{eq:22}) and part (iii) of
Theorem~\ref{pro:Main-rslt} in Lemma~\ref{sec:prel-notat} below. We
close this section by recalling some basic notations from Riemannian
geometry. In Section~\ref{sec:expansion-mu_2-small-5} we provide an
expansion of $\mu_2(B_g(y_0,r))$ as $r \to 0$. In
Section~\ref{sec:expansion-mu_2-small-3} we calculate a
corresponding expansion for suitably chosen geodesic ellipsoids with small
eccentricity. As shown by Corollary~\ref{sec:expansion-mu_2-small-2}, these
ellipsoids are suitable test domains to derive Part (i) of 
Theorem~\ref{pro:Main-rslt}, and from this the lower bound in Part
(ii) follows. Section~\ref{s:Ubm2} is devoted to collect all tools needed for
the proof of the upper bound in Theorem~\ref{pro:Main-rslt}(ii). In
particular, we use the above-mentioned stability estimate of Brasco and
Pratelli \cite{BP} in this
section, see Lemma~\ref{lem:asymd}.  Arguing by contradiction, we then
complete the proof of Theorem~\ref{pro:Main-rslt} in Section~\ref{sec:proof-main-result}.

\noindent
\bigskip

\textbf{Acknowledgments:}
This work is supported by the Alexander von Humboldt foundation. The
authors wish to thank Gerhard
 Huisken and the anonymous referee for helpful comments.

\section{Preliminaries and Notations}\label{s:pn}
We denote by  $B$ the unit ball in $\R^{N}$. Moreover, for a smooth
bounded domain $\Omega$ of a complete Riemannian manifold $(\M,g)$, we
write $\mu_2= \mu_2(\Omega,g)$ for the first nontrivial eigenvalue of
(\ref{eq:16}). The variational characterization of $\mu_2(\Omega,g)$ is given by 
\begin{equation}
  \label{eq:46}
\mu_2(\Omega,g) = \inf \biggl\{\frac{\int_{\Omega}|\nabla u|_g^2d v_g}{\int_{\Omega}
  u^2 dv_g}\::\: u \in H^1(\Omega) \setminus \{0\},\:
\int_{\Omega}u\,d v_g=0 \biggr\},   
\end{equation}
where $v_g$ denotes the volume element of the metric $g$. We recall
that the minimizers of this minimization problem are precisely the
eigenfunctions corresponding to $\mu_2(\Omega,g)$. If $\M= \R^{N}$ and $g$ is the euclidean metric, we
simply write $\mu_2(\Omega)$ in place of $\mu_2(\Omega,g)$.
As noted already, $\mu_2(B)$ is of multiplicity $N$ with corresponding eigenfunctions given by $\vp(|x|)\frac{x_i}{|x|}$, $i=1,\dots, N$ with
\be\label{eq:eqvp}
 \vp''+\frac{N-1}{t}\vp'+\left(\mu_2(B)-\frac{N-1}{t^2}
\right)\vp=0,\quad t\in (0,1),\quad \vp(0)=\vp'(1)=0.
 \ee
Throughout this paper, we assume the normalization \eqref{eq:normvp},
which equivalently yields 
\be\label{eq:normvpxiomx}
\int_B \vp^2(|x|)\,dx = N \qquad \text{and}\qquad  \int_B\vp^2(|x|)\left(\frac{x^i}{|x|}\right)^{2}dx={1} \quad
\text{for $i=1,\dots,N.$}
\ee
The function $ \vp$  and the eigenvalue $\mu_2(B)$ are obtained via
$J_{N/2}$, the Bessel function of the first kind of order $N/2$.
Indeed, $\sqrt{\mu_2(B)}$ is the first positive zero of the derivative
of  $t \mapsto t^{(2-N)/2} J_{N/2}( t)$, and $\vp$ is a scalar multiple of the function 
\begin{equation}
  \label{eq:42}
t \mapsto g(t)=t^{(2-N)/2} J_{N/2}(\sqrt{\mu_2(B)} t).
\end{equation}
More precisely, by \eqref{eq:normvp} we have
\begin{equation}
  \label{eq:23}
\vp(t)=\frac{g(t)}{\sqrt{|B|\int_0^1g^2 t^{N-1} dt}}.
\end{equation}

\begin{lem}
\label{sec:prel-notat}
We have:
\begin{itemize}
\item[(i)] $\mu_2(B) \ge N+1$ for every dimension $N \in \N$, and 
  \begin{equation}
    \label{eq:25}
\mu_2(B) < 
\left \{
\begin{aligned}
&N+2,&& \qquad \qquad \text{for $N=2,3,4$,}\\
&\frac{N(N-1)}{N-2} &&\qquad \qquad \text{for $N \ge 5$.}  
\end{aligned}
\right.    
\end{equation}
\item[(ii)] $|B| \vp^2(1)=
\frac{2\mu_2(B)}{\mu_2(B)-N+1}$ for every $N \in \N$.
\item[(iii)] $\gamma_N <0$ for all $N \ge 2$.
\item[(iv)] $\gamma_N \to 0$ as $N \to \infty$.
\end{itemize}
\end{lem}

\begin{pf}
Set $\mu_2:= \mu_2(B)$. We start with the proof of (ii). Since $\sqrt{\mu_2}$ is the first zero of the derivative of
the function $t \mapsto t^{(2-N)/2} J_{N/2}( t)$, we infer that
\begin{equation}
  \label{eq:39}
(J_{N/2}')^2(\sqrt{\mu_2})=\frac{(N-2)^2}{4}\frac{1}{\mu_2}J^2_{N/2}(\sqrt{\mu_2}).  
\end{equation}
Moreover, for the function $g$ defined in (\ref{eq:42}) we have by \cite[p.129, formula (5.14.5)]{lebedew}
\begin{equation}
  \label{eq:40}
 \int_0^1t^{N-1}g^2dt=\int_0^1t J^2_{N/2}(\sqrt{\mu_2}t)dt=\frac{1}{2}\Bigl[ 
 (J_{N/2}')^2(\sqrt{\mu_2})+\left( 1-\frac{N^2}{4 \mu_2} \right)J^2_{N/2}(\sqrt{\mu_2}) \Bigr].
\end{equation}
Inserting (\ref{eq:39}) in (\ref{eq:40}) yields
\begin{equation}
  \label{eq:43}
 \int_0^1t^{N-1}g^2dt= \frac{\mu_2-N+1}{2
   \mu_2}J^2_{N/2}(\sqrt{\mu_2}). 
\end{equation}
In particular, this implies $\mu_2 > N-1$. Moreover, since $g(1)=
J_{N/2}(\sqrt{\mu_2})$, we conclude by (\ref{eq:23}) that 
$$
|B|\vp^2(1)= \frac{g^2(1)}{\int_0^1t^{N-1}g^2 dt}= \frac{2 \mu_2}{\mu_2-N+1},
$$
as claimed.\\
We now turn to (i), and we first prove that $\mu_2 \ge N+1$. Let
$u,v: B \to \R$ be given by $u(x)=x_1$ and $v(x)= \frac{x_1}{|x|}\phi(|x|)$, so that
$-\Delta v = \mu_2 v$ in $B$ and $\partial_{\eta} v=0$ on $\partial
B$. Hence we find   
\begin{align*}
|B| \phi(1)&= \int_{\partial B}
\frac{x_1^2}{|x|}\phi(|x|)\,d\sigma = \int_{\partial B}uv\,d\sigma= \int_{\partial B}\Bigl(v \partial_{\nu} u- u \partial_{\nu}v\Bigr)\,d\sigma\\
&= \int_{B}\Bigl(v \Delta u- u \Delta v\Bigr)\,dx = \mu_2 \int_{B}u v\,dx=\frac{\mu_2}{N} \int_{B}|x|\phi(|x|)\,dx \\
&= \mu_2|B|\int_{0}^1 t^N
\phi(t)\,dt \le \mu_2 |B| \Bigl(\int_0^1 t^{N+1}\,dt
\Bigr)^{\frac{1}{2}} \Bigl(\int_0^1 t^{N-1}\phi^2(t) \,dt
\Bigr)^{\frac{1}{2}} = \frac{\mu_2 \sqrt{|B|}}{\sqrt{N+2}},  
\end{align*}
 using H{\"o}lder's inequality and \eqref{eq:normvp} in the last two
 steps. Using (ii) we therefore get 
$$
\mu_2^2-(N-1)\mu_2 -2(N+2) \ge 0,
$$
and this gives 
$$
\mu_2 \ge \frac{1}{2}\Bigl(N-1 + \sqrt{(N-1)^2+8(N+2)}\Bigr)\ge
 N+1.
$$ 
To prove (\ref{eq:25}), we consider the functions $x \mapsto
u_s(x)= x_1|x|^s$ for $s >-\frac{N}{2}$, so that $u \in H^1(B)$ and $\int_B
u_s\,dx = 0$. Since these functions are not eigenfunctions
corresponding $\mu_2$, the variational
characterization (\ref{eq:46}) gives
\begin{align*}
\mu_2 < \frac{\int_{B} |\nabla u_s|^2 \,dx}{\int_{B}
  u_s^2\,dx}=
\frac{\int_{B} \Bigl(|x|^{2s}+(s^2+2s)x_1^2|x|^{2(s-1)} 
\Bigr)\,dx}{
\int_B x_1^2|x|^{2s}\,dx}&=\frac{(N+s^2+2s
)\int_{B}|x|^{2s}\,dx}{\int_B |x|^{2s+2}\,dx}\\
&= \frac{(N+s^2+2s)(N+2s+2)}{N+2s}. 
\end{align*}
If $N \in \{2,3,4\}$, we may take $s=0$ and obtain the first inequality in
(\ref{eq:25}). If $N \ge 5$, we may take $s =-1$ and obtain the second
inequality in (\ref{eq:25}). This finishes the proof of (i).\\
To prove (iii), we consider the function 
$$
\sigma_N: (N-1,\infty) \to \R,\qquad \sigma_N(t)= \frac{1}{3N(N+2)} + \frac{N-2}{6Nt}
-\frac{1}{3N(t-N+1)}
$$
and recall from (\ref{eq:21}) that $\gamma_N= \sigma_N(\mu_2)$. We note
that 
\begin{equation}
  \label{eq:27}
\sigma_N'(t) = \frac{2t^2-(N-2)(t-N+1)^2}{6N(t-N+1)^2 t^2 }=
\frac{(4-N)t^2+(N-2)(N-1)[2t-(N-1)]}{6N(t-N+1)^2 t^2} .
\end{equation}
Hence $\sigma_N'(t) > 0$ in $(N-1,\infty)$ if $N \in \{2,3,4\}$, and therefore
(\ref{eq:25}) implies that 
$$
\gamma_N = \sigma_N(\mu_2) < \sigma_N(N+2) = 
 \frac{N-4}{18N(N+2)} \le 0 \qquad \text{if $N \in \{2,3,4\}$.}
$$
If $N \ge 5$, the zeros of the numerator in (\ref{eq:27}) are given by
$ \tau_N^\pm=\frac{N-1}{N-4}[N-2\,\pm\, \sqrt{2(N-2)}]$, whereas $\tau_N^- <
N-1 < \frac{N(N-1)}{N-2} <\tau_N^+$. Moreover, $\sigma_N'(t) > 0$ on 
$(N-1,\tau_N^+)$, and thus (\ref{eq:25}) implies
that  
\begin{align*}
\gamma_N&= \sigma(\mu_2)< \sigma_N(\frac{N(N-1)}{N-2}) 
= \frac{4-N}{3N^2(N+2)(N-1)}<0 \qquad \text{if $N \ge 5$.}
\end{align*}
This ends the proof of (iii).\\
(iv) simply follows from the definition of $\gamma_N$ and the fact
that $\gamma_N \ge N+1$, as shown in (ii).
\end{pf}

Let  $(\M,g)$ be a complete Riemannian manifold of dimension
$N$. We fix $y_0\in\M$ and consider an orthonormal basis
$E_1,\dots,E_{N}$ of $T_{y_0}\M$. In the sequel, it will be convenient to use the (somewhat sloppy) notation
$$
X:=x^i E_i \in T_{y_0}\M \qquad \text{for $x \in \R^{N}$.}
$$
Here and in the following, we sum over repeated upper and lower indices as usual. We
consider the geodesic coordinate system
\be\label{eq:sys-coor}
\R^{N}\ni x\mapsto\Psi(x):=\textrm{Exp}_{y_0}(X)
\ee
A geodesic ball in $\M$ centered at $y_0$ with radius $r>0$ is defined
as $B_g(y_0,r)=\Psi(rB)$. The map $\Psi$ induces coordinate
vector fields $Y_i:=\Psi_*  \frac{\partial}{\partial x^i}$, which
are pointwise given  by
$$
Y_i(x)= d\, \textrm{Exp}_{y_0}(X) E_i \in T_{\Psi(x)} \M, \qquad i=1,\dots,N.
$$
As usual, we write the metric in local coordinates by setting
$$
g_{ij}(x)=\la Y_i(x),Y_j(x)\ra_g \qquad \text{for $x \in \R^{N}$.}
$$
 The following local expansions are well known and can be found e.g. 
in \cite[\S II.8]{chavel-1}.
\begin{lem}  \label{l:oovg} In the above notations,
for any  $i,j=1,...,N$, we
have
$$
 g_{ij}(x)=\delta_{ij}+\frac{1}{3}\,\la R_{y_0}(X,E_i)X,E_j\ra_g
+ {O}(|x|^3)\quad \text{and}\quad
dv_{g}(x)=\Bigl( 1-\frac{1}{6}\, Ric_{y_0} (X,X) + O(|x|^3) \Bigr)dx.
$$
Here $dx$ is the volume element of $\R^{N}$,
$dv_g$ is the volume element of $\M$,
$$
R_{y_0}: T_{y_0}\M \times T_{y_0}\M \times T_{y_0}\M \to T_{y_0}\M
$$
is the Riemannian curvature tensor at $y_0$ and
$$
Ric_{y_0} : T_{y_0}\M \times T_{y_0}\M \to \R,\qquad Ric_{y_0}(X,Y)=-\sum_{i=1}^{N}\la R_{y_0}(X,E_i)Y,E_i\ra
$$
is the Ricci tensor at $y_0$. Moreover, the volume
expansion of metric balls is given by
\be\label{eq:expvolBgr}
\left|B_g(y_0,r) \right|_g= r^{N}\,{\left|B \right|}\,\left(1-\frac{1}{6(N+2)}\,r^2{S}(y_0)+O(r^4)  \right);
\ee
where ${S}$ is the scalar curvature function on $\M$.
\end{lem}
Here and in the following, once $y_0$ is fixed, we also write $\langle \cdot, \cdot
\rangle$ in place of $\langle \cdot, \cdot
\rangle_g$ to denote the scalar product on $T_{y_0}\M$ induced by the
metric $g$. It will turn out useful to put
\begin{equation}
  \label{eq:17}
R_{ijkl}:= \la R_{y_0}(E_i,E_j)E_k,E_l\ra  \quad \text{and}\quad R_{ij}:= Ric_{y_0}(E_i,E_j) \qquad \text{for $i,j=1,\dots,N$.}
\end{equation}
The scalar curvature of $\M$ at $y_0$ is given by $S(y_0)= \sum
\limits_{i=1}^{N}R_{ii}$. We point out that the orthonormal
basis $E_i$, $i=1,\dots,N$ can be chosen such that 
\begin{equation}
  \label{eq:37}
R_{ij}= 0  \qquad \text{for $i \not=j$},
\end{equation}
and we will fix such a choice from now on. We finally note that the euclidean
scalar product of $x,y \in \R^N$ will simply be denoted by $x \cdot
y$.

\section{Expansion of $\mu_2$ for small geodesic balls}
\label{sec:expansion-mu_2-small-5}
The main goal of the this section is the derivation of the the following expansion
for $\mu_2$ on small geodesic balls centered at $y_0$.

\begin{pro}\label{lem:expa_un2B}
For $r >0$ we have
$$
\mu_2(B_g(y_0,r),g) =\frac{\mu_2(B)}{r^2}+\a_N^-{S}(y_0) +2
\a_N^+R_{min}(y_0)+o(1),
$$
where
$$
\a_N^-=\frac{1}{6}\left(\frac{|B|\vp^2(1)}{N+2}-1\right),\quad
\a_N^+=\frac{1}{6}\left(\frac{|B|\vp^2(1)}{N+2}+1\right),\quad
R_{min}(y_0)=  \inf \limits_{A\in T_{y_0}{\M},|A|=1}
Ric_{y_0}(A,A)
$$
and $o(1) \to 0$ as $r \to 0$.
\end{pro}
\begin{pf}
Let $r>0$ be smaller than the injectivity radius of $\M$ at $y_0$, so
that $B_g(y_0,r)$ is a regular domain. Let $u_r \in C^3(\overline{B_g(y_0,r)})$ be an eigenfunction corresponding to the eigenvalue problem
$$
\D_g u_r+\mu_2(B_g(y_0,r),g) \,u_r=0\quad \text{in $B_g(y_0,r)$},\qquad  \la \n
u_r,\eta_r\ra_g=0\quad \text{on $\de B_g(y_0,r)$,}
$$
where $\eta_r$ denotes the outer unit normal on $\de B_g(y_0,r)$. Replacing $u_r$ by a scalar multiple if necessary, we may assume that
$u_r$ is a minimizer of the minimization problem
$$
\mu_2(B_g(y_0,r),g )=\inf\left\{ \int_{B_g(y_0,r) }\!\! |\n
  f|_g^2\,dv_g :
f \in H^1(B_g(y_0,r)), \int_{B_g(y_0,r) } \!\!f^2\,dv_g=1,\,\,  \int_{B_g(y_0,r) }\!\! f\,dv_g=0\right\}.
$$
Via the exponential map, we pull back the problem to the unit
ball $B \subset \R^{N}$. For this we consider the pull back metric
of $g$ under the map $B \to \M,\;x \mapsto \Psi(rx)$, rescaled with the
factor $\frac{1}{r^2}.$ Denoting this metric on $B$ by $g_r$, we then
have, in euclidean coordinates,
$$
[g_r]_{ij}(x)= \langle \frac{\de}{\de x^i},\frac{\de}{\de x^j}
\rangle_{g_r}\Big|_x  = \langle Y_i(\Psi(rx)),
Y_j(\Psi(rx))\rangle_{g}= g_{ij}(rx),
$$
so that
\begin{equation}
  \label{eq:1}
[g_r]_{ij}(x)=\d_{ij}+\frac{r^2}{3}\la R_{y_0}(X,E_i)X,E_j \ra+O(r^3)
\end{equation}
and
\begin{equation}
\label{eq:2}
g_r^{ij}(x)=\d^{ij}-\frac{r^2}{3}\la R_{y_0}(X,E_i)X,E_j \ra+O(r^3)
\end{equation}
uniformly for $x \in \overline B$ as a consequence of Lemma~\ref{l:oovg}. Here, as usual,
$(g_r^{ij})_{ij}$ denotes the inverse of the matrix
$([g_r]_{ij})_{ij}$. Setting $|g_r|= \det
([g_r]_{ij})_{ij}$, we also have $\sqrt{|g_r|}(x) = 1-\frac{r^2}{6}
Ric_{y_0}(X,X) +O(r^3)$ for $x \in \overline B$ by Lemma~\ref{l:oovg}. Since this expansion is valid
in the sense of $C^1$-functions on $\overline B$, we have
\begin{equation}
\frac{\de}{\de x^i}\sqrt{|g_r|}= -\frac{r^2}{3}Ric_{y_0}(X,E_i)+O(r^3) \qquad \text{for $i=1,\dots,N$.}  \label{eq:3}
\end{equation}
We now consider the rescaled eigenfunction
$$
\Phi_r: \overline B \to \R,\qquad \Phi_r(x)=r^{\frac{N}{2}}u_r(\Psi(rx))
$$
which satisfies 
$$
\D_{g_r} \Phi_r+\mu_2(B,g_r ) \,\Phi_r=0\quad \text{in $B$},\qquad
\la \n \Phi_r,\eta \ra_{g_r}=0\quad \text{on $\de B$,}
$$
with 
$$
\D_{g_r}\Phi_r=  \frac{1}{{\sqrt{|g_r|}}}\frac{\de }{\de x^i}\left(\sqrt{|g_r|}g_r^{ij} \frac{\de \Phi_r}{\de x^j}\right)
\qquad \text{and}\qquad \mu_2(B,g_r)=r^2 \mu_2(B_g(y_0,r),g).
$$
Moreover, $\int_{B} \Phi_r^2\,dv_{g_r}=1$ and $\int_{B}
\Phi_r\,dv_{g_r}=0$ with $ dv_{g_r}= \sqrt{|g_r|}dx$. Since
$g_r$ converges to the Euclidean metric in $\ov{B}$, it is easy to see
from the variational characterization of $\mu_2$ that
$\mu_2(B,g_r )\to \mu_2(B)$. Moreover, by using standard elliptic
regularity theory and compact Sobolev embeddings, one may show that,
along a sequence $r_k \to 0$, we have
$\Phi_{r_k}\to \Phi$ in $H^1(B)$ for some function $\Phi \in
C^2_{loc}(B)\cap C^1(\ov{B})$ satisfying
$$
\D \Phi+\mu_2(B ) \,\Phi=0\quad \text{in $B$} ,\qquad  \la \n \Phi,\eta
  \ra=0\quad \text{on $\de B$}, \quad \int_{B} \Phi^2\,dx=1 \quad
  \text{and}\quad  \int_{B} \Phi\,dx=0.
$$
Hence there exists $a=(a_1,\dots,a_N)=(a^1,\dots,a^N) \in \R^{N}$
with $|a|=1$ and such that
$$
\Phi(x)=\varphi(|x|)\frac{a  \cdot x }{|x|} \qquad \text{for $x
  \in \overline B$.}
$$
For matters of convenience, we will continue to write $r$ instead of
$r_k$ in the following. By integration by parts, using $\la \n\Phi_r,\eta\ra_{g_r}=0$ and
$dv_{g_r}= \sqrt{|g_r|}dx$, we have
$$
\mu_2(B, g_r)\, \int_{B} \Phi \Phi_r dv_{g_r}= -\int_{B} \Phi \D_{g_r}\Phi_r dv_{g_r}=
 \int_{B}\sqrt{|g_r|}g_r^{ij} \frac{\de \Phi}{\de x^i}
                                          \frac{\de \Phi_r}{\de x^j} \,dx.
$$
In the following, it will be convenient to use the notation
$$
\tn h = \sum_{i=1}^N \frac{\de h}{\de
  x^i}E_i: \overline B \to T_{y_0}\M \qquad \text{for a $C^1$-function
  $h$ defined on $\overline B$.}
$$
With this notation we find, using (\ref{eq:2}) and (\ref{eq:3}) and integrating by parts again,
\begin{align}
&\int_{B}\sqrt{|g_r|}g_r^{ij} \frac{\de \Phi}{\de x^i}
                                           \frac{\de \Phi_r}{\de x^j}
                                           \,dx
= \int_{B} \sqrt{|g_r|}\Bigl( \nabla \Phi_r \cdot \nabla \Phi -\frac{r^2}{3} \int_{B}\la
                                         R_{y_0}(X,\tn\Phi_r)X,\tn\Phi\ra\Bigr)dx+O(r^3)
\label{eq:4}\\
&=- \int_{B}\sqrt{|g_r|}(\D\Phi)\Phi_r\,dx
                                  - \int_{B}\Phi_r  \nabla
                                  \sqrt{|g_r|} \cdot \nabla \Phi  \,dx
                                          -\,\frac{r^2}{3} \int_{B}\la
                                          R_{y_0}(X,\tn\Phi_r)X,\tn\Phi\ra\,dx+O(r^3)
                                          \nonumber \\
&=\mu_2(B)\int_{B}\Phi\Phi_r\,dv_{g_r}+\frac{r^2}{3}\int_{B}\Phi_r\, Ric_{y_0}(X,\tn\Phi)\,dx-\,\frac{r^2}{3} \int_{B}\la R_{y_0}(X,\tn\Phi_r)X,\tn\Phi\ra\,dx+O(r^3).\nonumber
\end{align}
Therefore, since $ \int_{B}\Phi\Phi_r\,dv_{g_r}\to  1$ and $\Phi_r\to \Phi$ in $H^1(B)$ as $r\to0$, we obtain
\begin{equation}
  \label{eq:34}
\mu_2(B, g_r)=\mu_2(B)+\frac{r^2}{3}\int_{B}\Phi\, Ric_{y_0}(X,\tn\Phi)\,dx -\frac{r^2}{3} \int_{B}\la R_{y_0}(X,\tn\Phi)X,\tn\Phi\ra\,dx+o(r^2).
\end{equation}
Noticing that
\be\label{eq:nvp}
\tn \Phi(x)= \frac{1}{|x|^2} \Bigl(\vp'(|x|)-\frac{\vp(|x|)}{|x|}\Bigr) \la A,X \ra X
+\frac{\vp(|x|)}{|x|}A \qquad\quad \text{with $A := a^i E_i \in T_{y_0}\M$},
\ee
we find 
\begin{align}
\int_{B}\la R_{y_0}(X,\tn\Phi)X,&\tn\Phi\ra\,
dx=\int_{B}\frac{\vp^2(|x|)}{|x|^2} \la
R_{y_0}(X,A)X,A\ra \, dx=\la
R_{y_0}(E_i,A)E_j,A \ra  \int_{B}\frac{\vp^2(|x|)}{|x|^2}x^i x^j\,dx \nonumber\\ 
&=\la R_{y_0}(E_i,A)E_j,A \ra \int_{0}^1{\vp^2}{t^{N-1}}\,dt \int_{\partial B}x^i x^j \,d\s= -Ric_{y_0}(A,A).
 \label{eq:36}
\end{align}
Here we used the identity $\int_{\partial B}x^i x^j\,d\s =\delta^{ij}|B|$
and the normalization \eqref{eq:normvp} in the last step. Moreover, we compute via integration by parts, using (\ref{eq:37}), 
\begin{align}
2\int_{B}\Phi \,Ric_{y_0}(X,\tn\Phi)\,dx&= 2 R_{ij}\int_{B}\! x^i
\frac{\de \Phi}{\de x^j} \Phi\,dx=R_{ii}\int_{B}
x^i \frac{\de \Phi^2}{\de x^i}\,dx \nonumber\\
&= -\Bigl(\sum_{i=1}^{N} R_{ii}\Bigr)
\int_B \Phi^2\,dx + R_{ii} \int_{\de B}[x^i]^2 \Phi^2\,d\s \nonumber\\
&= -S(y_0)
\int_B \frac{(a \cdot x)^2}{|x|^2} \phi^2(|x|)\,dx + \phi^2(1) R_{ii} \int_{\de
  B}[x^i]^2 (a \cdot x)^2 \,d \sigma. \nonumber
\end{align}
Recalling \eqref{eq:normvpxiomx} and using the identities 
$$
\int_{\partial B}\![x^1]^4d\s =3\!\! \int_{\partial B}\! [x^i]^2 [x^j]^2 d\s =
\frac{3|B|}{N+2} \;\, \text{and}\; \int_{\partial B}\!x^i x^j
[x^k]^2 d\s=0\quad \text{for $i,j,k=1,\dots,N$, $i \not=j$,}
$$
we find that
\begin{align}
2&\int_{B}\Phi \,Ric_{y_0}(X,\tn\Phi)\,dx= -S(y_0) + \phi^2(1)
R_{ii}a_k a_l \int_{\de
  B}[x^i]^2  x^k  x^l \,d \sigma \nonumber\\
&= -S(y_0) + \vp^2(1)R_{ii}[a_k]^2 \int_{\de
  B}[x^i]^2  [x^k]^2\,d \sigma \nonumber= -S(y_0) + \frac{\vp^2(1)|B|}{N+2}\Bigl(S(y_0)+ 2
  R_{kk} [a^k]^2 \Bigr) \nonumber\\   
&=\Big(\frac{\vp^2(1)|B|}{N+2}-1\Big)\,S(y_0)+2 \frac{\vp^2(1)|B|}{N+2}\, Ric_{y_0}(A,A).\label{eq:14}
\end{align}
Combining (\ref{eq:34}), (\ref{eq:nvp}) and (\ref{eq:36}), we get
$$
\mu_2(B, g_r)=\mu_2(B) +\frac{r^2}{6}\Big(\frac{|B|\vp^2(1)}{N+2}-1
\Big)S(y_0)+\frac{r^2}{3}\Big(\frac{|B|\vp^2(1)}{N+2}+1
\Big)Ric_{y_0}(A,A) + o(r^2)
$$
and therefore
\begin{equation}
  \label{eq:10}
\mu_2(B_g(y_0,r),g) =\frac{\mu_2(B,g_r)}{r^2}= \frac{\mu_2(B)}{r^2}+\a_N^-{S}(y_0) +2
\a_N^+Ric_{y_0}(A,A)  +o(1).
\end{equation}
We now need to recall that -- more precisely -- here we have passed to
a sequence $r=r_k \to 0$. Nevertheless, the argument implies that
\begin{equation}
  \label{eq:8}
\mu_2(B_g(y_0,r),g) \ge \frac{\mu_2(B)}{r^2}+\a_N^-{S}(y_0) +2
\a_N^+R_{min}(y_0)+o(1) \qquad \text{as $r \to 0$.}
\end{equation}
Indeed, if - arguing by contradiction - there is
a sequence $r_k \to 0$ such that
\begin{equation}
  \label{eq:9}
\limsup_{k \to \infty} \Bigl[\mu_2(B_g(y_0,r_k),g) -
\frac{\mu_2(B)}{r_k^2}\Bigr] < \a_N^-{S}(y_0) +2
\a_N^+R_{min}(y_0),
\end{equation}
then by the above argument there exists a subsequence along which the
expansion (\ref{eq:10}) holds with some $A \in T_{y_0}\M$ with $|A|=1$, thus
contradicting (\ref{eq:9}). By (\ref{eq:8}), the proof of
Proposition~\ref{lem:expa_un2B} is finished once we have shown that
\begin{equation}
  \label{eq:44}
\mu_2(B_g(y_0,r),g) \leq\frac{\mu_2(B)}{r^2}+\a_N^-{S}(y_0) +2
\a_N^+Ric_{y_0}(A,A) +o(1)
\end{equation}
for all $A \in T_{y_0}\M$ with $|A|=1$. So now consider 
$a=(a^1,\dots,a^N) \in \R^N$ arbitrary with
$|a|=1$, and let $A=a^i E_i \in
T_{y_0}\M$. We define
$$
\ti{\Phi}: \overline B \to \R,\qquad \ti{\Phi}(x)=\varphi(|x|)\frac{a
  \cdot x}{|x|}
$$
and
$$
c_r:= \frac{1}{|B|_{g_r}}\int_{B} \ti{\Phi}dv_{g_r}.
$$
Then, by Lemma~\ref{l:oovg},
$$
c_r = \Bigl(\frac{1}{|B|}+O(r^{2})\Bigr)\Bigl(
 \int_{B}\ti{\Phi}(x)[1-\frac{1}{6} Ric_{y_0}(X,X)]dx
 +O(r^3)\Bigr)=\Bigl(\frac{1}{|B|}+O(r^{2})\Bigr)O(r^3)=O(r^3),
$$
since the function $x \mapsto \ti{\Phi}(x)[1-\frac{1}{6} Ric_{y_0}(X,X)]$ is odd with
respect to reflection at the origin. Hence, using the
variational characterization of $\mu_2(B,g_r)$, we find that
$$
\mu_2(B,g_r)\leq\,\frac{  \displaystyle \int_{B } |\n (\ti{\Phi}-c_r)|_{g_r}^2\, d{v_{g_r}}    }{\displaystyle  \int_{B}
 (\ti{\Phi}-c_r)^2 \,d{v_{g_r}}}
=\,\frac{  \displaystyle \int_{B } |\n \ti{\Phi}|_{g_r}^2\, d{v_{g_r}}
}{\displaystyle  \int_{B} \ti{\Phi}^2 +O(r^3) \,d{v_{g_r}}}
=\,\frac{  \displaystyle \int_{B } |\n \ti{\Phi}|_{g_r}^2\,d{v_{g_r}}  +O(r^3)}{\displaystyle  \int_{B}
 \ti{\Phi}^2 \,d{v_{g_r}}}
$$
and therefore
$$
\mu_2(B,g_r) \int_{B }
  \ti{\Phi}^2\,d{v_{g_r}} \le
 \int_{B_g(y_0,r) } |\n \ti{\Phi}|_{g_r}^2\,d{v_{g_r}} +O(r^3)=
\int_{B}\sqrt{|g_r|}g_r^{ij} \frac{\de \ti {\Phi}}{\de x^i}
                                           \frac{\de \ti{\Phi}}{\de x^j}
                                           \,dx  +O(r^3).
$$
It is by now straightforward that the same estimates as above -- starting from
(\ref{eq:4}) -- hold with both $\Phi_r$ and $\Phi$ replaced by
$\ti{\Phi}$. We thus obtain (\ref{eq:44}), as required.  
\end{pf}

\begin{cor}
\label{sec:expansion-mu_2-small-1}
We have
\be\label{eq:expmu_2}
\mu_2(B_g(y_0,r),g)=\left(1- \b(y_0)\left( \frac{v}{|B|}
\right)^{\frac{2}{N}}   +o\left( \frac{v}{|B|}
\right)^{\frac{2}{N}} \right) SW_{\R^{N}}(v)
\ee
as $v=\left| B_g(y_0,r)\right|_g \to 0$ with
$$
\b(y_0) = \frac{\mu_2(B)S(y_0)  -3N(N+2)\left(
\a_N^-{S}(y_0) +2 \a_N^+R_{min}(y_0) \right) }{ 3N(N+2)\mu_2(B)
}.
$$
\end{cor}

\begin{pf}
By the volume expansion (\ref{eq:expvolBgr}) of geodesic balls we have
\begin{align*}
\Bigl(\frac{v}{r^{N}
  |B|}\Bigr)^{\frac{2}{N}}=\Bigl(\frac{|B_g(y_0,r)|_g}{r^{N}
  |B|}\Bigr)^{\frac{2}{N}}&= 1 - \frac{1}{3N(N+2)}S(y_0)r^2+o(r^2)\\
&= 1 - \frac{1}{3N(N+2)}S(y_0)\Bigl(\frac{v}{|B|}\Bigr)^{\frac{2}{N}} +o\Bigl(\frac{v}{|B|}\Bigr)^{\frac{2}{N}}   
\end{align*}
as $v=\left| B_g(y_0,r)\right|_g \to 0$. Together with Lemma~\ref{lem:expa_un2B} this yields 
\begin{align*}
\mu_2&(B_g(y_0,r),g) =\frac{\mu_2(B)}{r^2}+\a_N^-{S}(y_0) +2
\a_N^+R_{min}(y_0)+o(1)\\
&= \left( \left(\frac{v}{r^{N} |B|}\right)^{\frac{2}{N}} +\frac{\a_N^-{S}(y_0) +2
\a_N^+R_{min}(y_0)}{\mu_2(B)}\left( \frac{v}{|B|}
\right)^{\frac{2}{N}} +o\left( \frac{v}{|B|}
\right)^{\frac{2}{N}}  \right) SW_{\R^{N}}(v)\\
&= \left(1-\left(\frac{1}{3N(N+2)} S(y_0)- \frac{\a_N^-{S}(y_0) +2
\a_N^+R_{min}(y_0)}{\mu_2(B)}\right)\left( \frac{v}{|B|}
\right)^{\frac{2}{N}} +o\left( \frac{v}{|B|}
\right)^{\frac{2}{N}}  \right) SW_{\R^{N}}(v)\\
&=\left(1- \b(y_0)\left( \frac{v}{|B|}
\right)^{\frac{2}{N}}   +o\left( \frac{v}{|B|}
\right)^{\frac{2}{N}} \right) SW_{\R^{N}}(v)
\end{align*}
as $v=\left| B_g(y_0,r)\right|_g \to 0$.
\end{pf}

\begin{rem}
\label{sec:expansion-mu_2-small}
{\rm (i) Since
\be\label{eq:RicleS}
R_{min}(y_0)\leq \frac{S(y_0)}{N},
\ee
and 
\begin{equation}
  \label{eq:41}
\alpha_N^- + \frac{2 \alpha_N^+}{N}= \frac{|B| \vp^2(1)-(N-2)}{6N},  
\end{equation}
Proposition~\ref{lem:expa_un2B} and Corollary~\ref{sec:expansion-mu_2-small-1} yield 
\begin{align}
\mu_2(B_g(y_0,r),g) &\le \frac{\mu_2(B)}{r^2}+(\a_N^-+
\frac{2\a_N^+}{N}){S}(y_0)+o(1)   \label{eq:45}\\
&= \left(1-\g_{N}\left(\frac{v}{|B|}\right)^{\frac{2}{N}}\,S(y_0)
 +o\left( \frac{v}{|B|}
\right)^{\frac{2}{N}}    \right)\,SW_{\R^{N}}(v)\nonumber
\end{align}
as $v=\left| B_g(y_0,r)\right|_g \to 0$ (and therefore $r \to 0$) with
$\gamma_N$ as in (\ref{eq:21}). Notice that when $N=2$, equality holds in \eqref{eq:RicleS} and (\ref{eq:45}). Therefore the two-dimensional version  of \eqref{eq:expmu_2} is
$$
\mu_2(B_g(y_0,r),g)=\left(1- \g_2 \frac{v}{|B|}
\, S(y_0)  +o\left( \frac{v}{|B|}
\right)  \right) SW_{\R^{2}}(v). 
$$
(ii) Denote by $(\mathbb{M}^{N},g_k)$ a space of constant
sectional curvature $k$. Then equality holds in \eqref{eq:RicleS} because $Ric=(N-1)k\,g_k$ on $\mathbb{M}^{N}$. In
particular if $E$ is a ball in $(\mathbb{M}^{N},g_k)$ with small
volume, one has that
\be\label{eq:mu2Ek}
 \mu_2(E,g_k)= \left(1-\g_{N}\left(\frac{|E|_{g_k}}{|B|}\right)^{\frac{2}{N}}\,N(N-1)\,k
 +o \left(\frac{|E|_{g_k}}{|B|}\right)^{\frac{2}{N}}  \right)\,SW_{\R^{N}}(| E|_{g_k}).
\ee
}
\end{rem}

\section{Expansion of $\mu_2$ for small geodesic ellipsoids}
\label{sec:expansion-mu_2-small-3}
As before we fix $y_0 \in \M$, and we continue to assume that the orthonormal basis
$E_1,\dots,E_N$ of $T_{y_0}\M$ is chosen such that (\ref{eq:37})
holds. In the following, we consider 
$$
\nu_N=\frac{2 \mu_2(B)+N|B|\vp^2(1)}{N+2}>0,
$$
and we let 
\begin{equation}
  \label{eq:30}
b_i=b^i:= \frac{\alpha_N^+}{\nu_N} (R_{ii} -\frac{S(y_0)}{N}) \qquad
\text{for $i=1,\dots,N$,}
\end{equation}
where $\alpha_N^+$ is defined in Proposition~\ref{lem:expa_un2B}. 
The reason for this choice will become clear later. We note that
$\sum \limits_{i=1}^N b_i= 0$ since $S(y_0)= \sum \limits_{i=1}^N
R_{ii}$. For $r>0$, we now consider the geodesic ellipsoids
$E(y_0,r):= F_r(B) \subset \M$, where 
$$
F_r: B \to \M, \qquad F_r(x)= \textrm{Exp}_{y_0}\bigl(r (1+r^2b_i)x^i E_i\bigr).
$$
The special choice of the values $b_i$ gives rise to the following
asymptotic expansion where the local geometry only enters via the
scalar curvature at $y_0$.

\begin{pro}
\label{sec:expansion-mu_2-small-4}
As $r \to 0$, we have
\begin{equation}
  \label{eq:29}
\mu_2(E(y_0,r),g)
=\frac{\mu_2(B)}{r^2}+(\alpha_N^- + \frac{2\alpha_N^+ }{N} )S(y_0)+o(1),
\end{equation}
with $\alpha_N^\pm$ as in Proposition~\ref{lem:expa_un2B} and 
\begin{equation}
  \label{eq:31}
|E(y_0,r)|_g=|B_g(y_0,r)|_g+O(r^{N+4})= r^N |B| \Bigl(1 -
\frac{1}{6(N+2)}r^2 S(y_0) + O(r^4)\Bigr).
\end{equation}
\end{pro}

\begin{pf}
We consider the pull back metric $h_r$ on $B$ 
of $g$ under the map $F_r$ rescaled with the
factor $\frac{1}{r^2}.$ Then we have 
\begin{align}
\label{eq:exp-hr}
[h_r]_{i,j}(x)&=(1+r^2 b_i)(1+r^2 b_j)[g_r]_{ij}((1+r^2 b_k)x^k
e_k)=[g_r]_{ij}(x)+ r^2(b_i+b_j)\d_{ij}+O(r^4)\\
&= \delta_{ij} + r^2 \Bigl(\frac{1}{3}\la R_{y_0}(X,E_i)X,E_j \ra + (b_i+b_j)\d_{ij}\Bigr)+O(r^3)\nonumber
\end{align}
uniformly in $x \in B$. Setting $|h_r|= \det ([h_r]_{ij})_{ij}$, we deduce the expansion 
$$
|h_r|(x)= |g_r|(x) + 2 r^2 \sum_{i=1}^N b_i + O(r^4)= |g_r|(x)+ O(r^4)
\qquad \text{for $x \in B$.}
$$
This implies that 
$$
|E(y_0,r)|_g = r^N |B|_{h_r}= r^N \Bigl(|B|_{g_r}+O(r^4)\Bigr)
=|B_g(y_0,r)|_g+O(r^{N+4}),
$$
as claimed in (\ref{eq:31}). \\
We now turn  to \eqref{eq:29}. We first note
that $\mu_2(B,h_r)=r^2 \mu_2(E(y_0,r),g)$; therefore (\ref{eq:29})
is equivalent to 
\begin{equation}
  \label{eq:33}
\mu_2(B,h_r)
=\mu_2(B)+r^2 (\alpha_N^- + \frac{2\alpha_N^+ }{N} )S(y_0)+o(r^2).
\end{equation}
Let $\Phi_r$ be an eigenfunction for $\mu_2(B, h_r)$, normalized such
that $\int_{B}\Phi_r^2\,dv_{h_r}=1$ with $dv_{h_r}= \sqrt{|h_r|}dx$. Then we have
 $$
\D_{h_r} \Phi_r+\mu_2(B,h_r ) \,\Phi_r=0\quad \text{in $B$},\qquad
\la \n \Phi_r,\eta \ra_{h_r}=0\quad \text{on $\de B$,}
$$
where 
$$
\D_{h_r}\Phi_r=  \frac{1}{{\sqrt{|h_r|}}}\frac{\de }{\de x^i}\left(\sqrt{|h_r|}h_r^{ij} \frac{\de \Phi_r}{\de x^j}\right).
$$
Since
$h_r$ converges to the Euclidean metric in $B$, the variational
characterization of $\mu_2$ implies that
$\mu_2(B,h_r )\to \mu_2(B)$. Moreover, as in the proof of
Proposition~\ref{lem:expa_un2B} we have $\Phi_{r_k}\to \Phi$ in
$H^1(B)$ along a sequence $r_k \to 0$ with some function $\Phi \in
C^2_{loc}(B)\cap C^1(\ov{B})$ satisfying
$$
\D \Phi+\mu_2(B ) \,\Phi=0\quad \text{in $B$} ,\qquad  \la \n \Phi,\eta
  \ra=0\quad \text{on $\de B$}, \quad \int_{B} \Phi^2\,dx=1 \quad
  \text{and}\quad  \int_{B} \Phi\,dx=0.
$$
Hence there exists a vector $a=(a_1,\dots,a_{N})=(a^1,\dots,a^{N}) \in \R^{N}$
with $|a|=1$ and such that
$$
\Phi(x)=\varphi(|x|)\frac{a \cdot x}{|x|}   \qquad \text{for $x
  \in \overline B$.}
$$
For matters of convenience, we will continue to write $r$ instead of
$r_k$ in the following. By multiple integration by parts, using
\eqref{eq:exp-hr} and (\ref{eq:3}), we have 
\begin{align}
\mu_2(B, h_r)\, &\int_{B} \Phi \Phi_r dv_{h_r}= -\int_{B} \Phi \D_{h_r}\Phi_r dv_{h_r}=
 \int_{B}\sqrt{|h_r|}h_r^{ij} \frac{\de \Phi}{\de x^i}
                                          \frac{\de \Phi_r}{\de x^j}
                                          \,dx \nonumber \\
=& \int_{B} \sqrt{|h_r|}\Bigl[\nabla \Phi_r \nabla \Phi -r^2 \Bigl(\frac{1}{3}\la
                                         R_{y_0}(X,\tn\Phi_r)X,\tn \Phi\ra 
+ 2 b^i \frac{\de \Phi}{\de x^i}
                                           \frac{\de \Phi_r}{\de
                                             x^i}\Bigr)\Bigr]\,dx+O(r^3)\nonumber
                                           \\
=&- \int_{B}\sqrt{|h_r|}(\D\Phi)\Phi_r\,dx
                                  - \int_{B}\Phi_r \nabla \sqrt{|h_r|}
                                  \cdot \nabla \Phi \,dx -\,\frac{r^2}{3} \int_{B}\la
                                          R_{y_0}(X,\tn\Phi_r)X,\tn\Phi\ra\,dx
 \nonumber\\
&- 2r^2  \int_{B}b^i \frac{\de \Phi}{\de x^i}
                                           \frac{\de \Phi_r}{\de x^i}\,dx
                                          +O(r^3)
                                          \nonumber \\
=&\:
\mu_2(B)\int_{B}\Phi\Phi_r\,dv_{h_r}+\frac{r^2}{3}\int_{B}Ric_{y_0}(\tn\Phi,X)\Phi_r\,dx-
\frac{r^2}{3} \int_{B}\la R_{y_0}(X,\tn\Phi_r)X,\tn\Phi\ra\,dx
  \nonumber\\
&-2r^2  \int_{B}b^i \frac{\de \Phi}{\de x^i}
                                           \frac{\de \Phi_r}{\de x^i}\,dx+O(r^3).\nonumber
\end{align}
Since $ \int_{B}\Phi\Phi_r\,dv_{h_r}\to  1$ and $\Phi_r\to \Phi$ in
$H^1(B)$ as $r\to0$, we may use the calculations in the proof of Proposition
\ref{lem:expa_un2B} starting from (\ref{eq:34}) to obtain
\begin{align}
\mu_2(B, h_r)=&\:\mu_2(B)+\frac{r^2}{3}\int_{B}Ric_{y_0}(\tn\Phi,X)\Phi_r\,dx-
\frac{r^2}{3} \int_{B}\la
R_{y_0}(X,\tn\Phi_r)X,\tn\Phi\ra\,dx\nonumber \\
 &-2r^2  \int_{B}b^i \frac{\de \Phi}{\de x^i}
                                           \frac{\de \Phi_r}{\de
                                             x^i}\,dx+o(r^2)\nonumber \\ 
=&\: {\mu_2(B)}+{r^2}\Bigl( \a_N^-{S}(y_0) + 2
{r^2} \a_N^+Ric_{y_0}(A,A) - 2 \int_{B} b^i  \Bigl(\frac{\de
  \Phi}{\de x^i}\Bigr)^2\,dx \Bigr) +o(r^2)
\label{eq:32}
\end{align}
with $A := a^i E_i \in T_{y_0}\M$. It remains to compute $\int_{B} b^i  \Bigl(\frac{\de \Phi}{\de
  x^i}\Bigr)^2\,dx$.\\ We have $\frac{\partial \Phi}{\partial x^i}= \Bigl(\phi'(|x|)
-\frac{\phi(|x|)}{|x|}\Bigr)\frac{x^i}{|x|^2} a_j  x^j + a_i
\frac{\phi(|x|)}{|x|}$ and thus  
$$
\Bigl(\frac{\partial \Phi}{\partial x^i}\Bigr)^2 = \frac{1}{|x|^4}\Bigl(\phi'(|x|) - \frac{\phi(|x|)}{|x|} \Bigr)^2 a_j
a_k [x^i]^2 x^j x^k  + a_i^2 \frac{\phi^2(|x|)}{|x|^2} +2 \frac{\phi(|x|)}{|x|^3}
\Bigl(\phi'(|x|) - \frac{\phi(|x|)}{|x|}\Bigr) a_i a_j x^i x^j
$$
for $x \in B$ and $i=1,\dots,N$. Noting the oddness of some of the integrands and passing to polar
coordinates, we therefore obtain
\begin{align}
\int_{B} b^i   \Bigl(\frac{\de \Phi}{\de x^i}\Bigr)^2
\,dx&= b_i a_j^2 \int_{B} \frac{1}{|x|^4}\Bigl(
\phi'(|x|) - \frac{\phi(|x|)}{|x|} \Bigr)^2 
[x^i]^2 [x^j]^2\,dx+  b^i a_i^2 \int_{B} \frac{\phi^2(|x|)}{|x|^2}\,dx  \nonumber\\
&+ 2 b^i a_i^2 \int_B \frac{\phi(|x|)}{|x|^3}
\Bigl(\phi'(|x|) - \frac{\phi(|x|)}{|x|}\Bigr) [x^i]^2\nonumber
\\
&= b_i a_j^2 \int_0^1 t^{N-1} \Bigl(\phi'(t) - \frac{\phi(t)}{t} \Bigr)^2\;dt 
\int_{\partial B} [x^i]^2 [x^j]^2\,d\sigma + b^i a_i^2  |\partial B|
\int_{0}^1 t^{N-3} \phi^2\,dt\nonumber\\
&+ 2 b^i a_i^2 \int_0^1 t^{N-2} \phi(t) \Bigl(\phi'(t) -
\frac{\phi(t)}{t}\Bigr)\,dt \int_{\partial B}  [x^i]^2\,d \sigma. \label{eq:38}
\end{align}
Put $d_N:= \int_{0}^1 t^{N-3} \phi^2\,dt$. By \eqref{eq:eqvp}, we have
$$
\int_0^1t^{N-1} (\vp')^2(t) \,dt =\frac{\mu_2(B)}{|B|}-(N-1) d_N \qquad
\text{and}\qquad \int_0^1 t^{N-2} \vp'(t)\vp(t) \,dt =\frac{\vp^2(1)}{2}-\frac{N-2}{2} d_N
$$
hence 
$$
\int_0^1 t^{N-1} \Bigl(\phi'(t) - \frac{\phi(t)}{t} \Bigr)^2\;dt  =
\frac{\mu_2(B)}{|B|}- \vp^2(1) \qquad \text{and}\qquad \int_0^1 t^{N-2} \phi(t)\Bigl(\phi'(t) -
\frac{\phi(t)}{t}\Bigr)\,dt = \frac{1}{2}\Bigl(\vp^2(1) -N d_N\Bigr).
$$
Inserting this in (\ref{eq:38}), we get
\begin{align*}
\int_{B} b^i   \Bigl(\frac{\de \Phi}{\de x^i}\Bigr)^2
\,dx&= b_i a_j^2 \Bigl(\frac{\mu_2(B)}{|B|}- \vp^2(1)\Bigr)
\int_{\partial B} [x^i]^2 [x^j]^2\,d\sigma 
+ b^i a_i^2 |\partial B| d_N + b^i a_i^2\Bigl(\vp^2(1)-Nd_N\Bigr) 
 \int_{\partial B}  [x^i]^2\,d \sigma.
\end{align*}
Recalling furthermore the identities 
\begin{align*}
\int_{\de B} [x^i]^4 \,d\sigma = 3 \int_{\de B} [x^i]^2 [x^j]^2 \,d\sigma =
\frac{3}{N+2}|B|,\qquad \int_{\de B} (x^i)^2\,d\sigma=\frac{|\de B|}{N}=|B|
\end{align*}
for $i,j=1,\dots,N$, $i \not = j$ and also that $\sum \limits_{i=1}^N
b_i = 0$, we obtain 
$$
b_i \int_{\partial B} [x^i]^2 [x^j]^2\,d\sigma =
\frac{2|B|}{N+2}b^j\qquad \text{for $j=1,\dots,N$}
$$
and thus
\begin{align*}
\int_{B} b^i   \Bigl(\frac{\de \Phi}{\de x^i}\Bigr)^2
\,dx&= b^j a_j^2 \Bigl(\frac{\mu_2(B)}{|B|}- \vp^2(1)\Bigr) \frac{2|B|
}{N+2}
+ b^i a_i^2 N|B| d_N + b^i a_i^2\Bigl(\vp^2(1)
-N d_N\Bigr)|B|\\
&= b^i a_i^2 \,\frac{2 \mu_2(B)+N|B|\vp^2(1)}{N+2}= \nu_N\, b^i a_i^2 .
\end{align*}
Inserting this in (\ref{eq:32}), we obtain
\begin{align*}
\mu_2&(B, h_r)
=\: {\mu_2(B)}+{r^2}\Bigl[ \a_N^-{S}(y_0) + 2(
\a_N^+Ric_{y_0}(A,A) - \nu_N b^i a_i^2) \Bigr] +o(r^2)\\
=&\: {\mu_2(B)}+{r^2}\Bigl[ (\a_N^-+ \frac{2 a_N^+}{N}){S}(y_0) +
2(\a_N^+ (Ric_{y_0}(A,A)- \frac{{S}(y_0)}{N}) - \nu_N b^i a_i^2 )  \Bigr] +o(r^2),
\end{align*}
where 
$$
\a_N^+ (Ric_{y_0}(A,A)- \frac{{S}(y_0)}{N}) - \nu_N b^i a_i^2 =  [a^i]^2
\Bigl(\a_N^+ (R_{ii}- \frac{{S}(y_0)}{N}) - \nu_N b_i\Bigr)=0
$$
by our choice of the $b_i=b^i$ in (\ref{eq:30}). This shows
(\ref{eq:33}), as required.
\end{pf}

\begin{cor}
\label{sec:expansion-mu_2-small-2}
We have
\be\label{eq:expmu_2-ellipsoid}
\mu_2(E(y_0,r),g)=\left(1- \gamma_N \left( \frac{v}{|B|}
\right)^{\frac{2}{N}}   +o\left( \frac{v}{|B|}
\right)^{\frac{2}{N}} \right) SW_{\R^{N}}(v),
\ee
as $v=\left| E(y_0,r)\right|_g \to 0$ with $\gamma_N$ as in (\ref{eq:21}).
\end{cor}

\begin{pf}
This follows readily by combining (\ref{eq:41}), (\ref{eq:29}) and (\ref{eq:31}). 
\end{pf}

%
\section{A local upper bound for $\mu_2$}\label{s:Ubm2}
We fix $r_0>0$ less than the convexity radius of
$\M$ at $y_0$, so that $r_0$ is also less than the injectivity radius
of $\M$ at $y_0$. As in
\cite{Wein}, we consider the function
$$
G: \R \to \R,\qquad  G(t)= \left\{
\begin{aligned}
\vp(t)\quad\textrm{ if } t\leq 1,\\
\vp(1)\quad\textrm{ if } t>1,
\end{aligned}
\right.
$$
where $\vp$ is the function defined in Section~\ref{s:pn}. Throughout this section, we consider a sequence of numbers $r_k \in
(0,\frac{r_0}{3})$ such that $r_k \to 0$ as $k \to \infty$, and we
suppose that we are given regular domains $\O_{r_k}\subset
B_g(y_0,r_k)$, $k \in \N$. In order to keep the notation as simple as possible, we
will write $r$ instead of $r_k$ in the following. By \cite[Theorem 3]{AS}, there exists a point
$p_r\in B_g(y_0,r)$ such that
\be\label{eq:centmass}
\int_{\O_r}\frac{G(|\textrm{Exp}_{p_r}^{-1}(q)|_g)}{|\textrm{Exp}_{p_r}^{-1}(q)|_g  }\textrm{Exp}_{p_r}^{-1}(q)\,dv_g=0.
\ee
Moreover, there exists a unique $\rho_r \in (0,r)$ such that
that $|\O_r|_g=|B_g(p_r,\rho_r) |_g$. We have that, for every $r>0$ small,
 $B_g(p_r,\rho_r)\subset B_g(y_0,2r)$ and also $\O_r\subset B_g(p_r,2 r)$.
Now we need to extend some of the notations introduced
in Section~\ref{s:pn}. For this we let
$$
y \mapsto E_i^y \in T_y \M,\qquad i=1,\dots,N
$$
denote a smooth orthonormal frame on $B_g(y_0,r_0)$, and we define
$$
\Psi_r: \R^{N} \to \M,\qquad \Psi_r(x)= \textrm{Exp}_{p_r}(x^i E_i^y).
$$
We also define
 \be\label{eq:defUr}
B^r:= \frac{2r}{\rho_r}B
\qquad \text{and}\qquad
U_r:=\frac{1}{\rho_r}\,\Psi_r^{-1} (\O_r)
\subset B^r,
\ee
and we consider the
pull back metric of $g$ under the map $B^r \to \M,\; x \mapsto
\Psi_r(\rho_r x)$, rescaled with the factor $\frac{1}{\rho_r^2}$.
We denote this metric on $B^r$ by $g_r$, and we point out that this
definition differs from the notation used in the proof of
Proposition~\ref{lem:expa_un2B}. By \eqref{eq:centmass}, it is plain that
\be\label{eq:choOr}
\int_{U_r}\frac{G(|x|)}{|x|}x^idv_{g_r}=0 \qquad \text{for $i=1,\dots, N$.}
\ee
We also write
$$
R_{ijkl}^{r}:= \la R_y(E_i^{p_r},E_j^{p_r})E_k^{p_r},E_l^{p_r} \ra \quad \text{and}\quad
R_{ij}^{r}:= Ric_{{p_r}}(E_i^{p_r},E_j^{p_r})
$$
for $i,j,k,l=1,\dots,N$. To be consistent with the notation
introduced in the end of Section~\ref{s:pn},  we also write
$$
R_{ijkl}:= \la R_{y_0}(E_i^{y_0},E_j^{y_0})E_k^{y_0},E_l^{y_0} \ra \quad \text{and}\quad
R_{ij}:= Ric_{y_0}(E_i^{y_0},E_j^{y_0}). 
$$
Since $\dist(p_r,y_0)=O(r)$, we then have
\begin{equation}
  \label{eq:11}
R_{ijkl}^{r} =R_{ijkl}+O(r) \quad \text{and}\quad
R_{ij}^{r}=R_{ij}+O(r)  \qquad \text{for $i,j,k,l=1,\dots,N$.}
\end{equation}
By Lemma~\ref{l:oovg} we also have
\be\label{eq:metexp}
\begin{array}{rllll}
\displaystyle (g_r)_{ij}(x)&=\delta_{ij}+\frac{\rho_r^2}{3}
R_{kilj}^{r}x^kx^l +|x|^3 {O}(\rho_r^3);\\[3mm]
 \displaystyle dv_{g_r}(x)=\sqrt{|g_r(x)|}\,dx&=\left(
   1-\frac{\rho_r^2}{6}\, R_{lk}^{r}x^lx^k +|x|^3 {O}(\rho_r^3) \right)dx,
\end{array}
\ee
uniformly on $B^r$, where $|g_r|$ is the
determinant of $g_r$, so in particular
\be
\label{eq:metexp-1}
(g_r)_{ij}(x)=\delta_{ij}+O(r^2) \quad \text{and}\quad
dv_{g_r}(x)=(1+O(r^2))dx \qquad \text{uniformly on $B^r$.}
\ee
Observe that 
\be\label{eq:volU}
|U_r|_{g_r}={\rho_r}^{-N}|\Omega_r|_{g}={\rho_r}^{-N}|B_g(p_r,\rho_r) |_g=
|B|_{g_r} \qquad \text{and}\qquad \mu_2(U_r,g_r)=\frac{\mu_2(\O_r,\,g)}{\rho_r^2}. 
 \ee
Moreover, since $U_r \subset B^r$ and $B \subset B^r$, we infer from \eqref{eq:metexp-1} that
\be\label{eq:volUdx}
|U_r|=(1+O(r^2))|U_r|_{g_r}=(1+O(r^2))|B|_{g_r}=(1+O(r^2))|B|.
\ee
Setting $$
{ f}_i: \R^{N} \to \R,\qquad f_i(x)=\frac{G(|x|)}{|x|}x^i,
$$
we find that $\int_{U_r}f_idv_{g_r}=0$ for $i=1,\dots,N$ by \eqref{eq:choOr},
and hence the variational characterization of $\mu_2$ yields
\be\label{eq:mudlesSm}
  \mu_2(U_r,g_r)\leq\frac{\displaystyle
\sum_{i=1}^N\int_{U_r}|\n f_i|_{g_r}^2dv_{g_r}  }{\displaystyle
\sum_{i=1}^N \int_{U_r}f_i^2dv_{g_r}}.
 \ee
We also note that
\begin{equation}
  \label{eq:18}
 \frac{\de f_i}{\de x^k} =
\frac{G'}{|x|^2}x^ix^k   +\frac{G}{|x|}[\delta_{ik}-\frac{x^ix^k}{|x|^2}]
\qquad \text{for every $i,k = 1,\dots,N$}
\end{equation}
and, by direct calculation as in \cite{Wein},
\be\label{smfi2}
\sum_{i=1}^{N} f_i^2 =G^2,  \qquad \sum_{i=1}^{N}
|\n f_i|^2
 = (G')^2+(N-1) \frac{G^2}{|x|^2}.
\ee
Here and in the following, we simply write $G$
instead of $G(|\cdot|)$ or $G(|x|)$ and $G'$ instead of $G'(|\cdot|)$
or $G'(|x|)$ if the meaning is clear from the context.
In particular, using \eqref{smfi2}, \eqref{eq:normvpxiomx} and recalling that $\vp$ and $G$
coincide in $[0,1]$, we observe that
\begin{align}
  \label{eq:19}
\int_B \Bigl((G')^2+(N-1) \frac{G^2}{|x|^2}\Bigr)\,dx = \sum_{i=1}^{N} \int_B
|\n f_i|^2\,dx &= \mu_2(B) \sum_{i=1}^{N} \int_B
|f_i|^2\,dx\\
&= \mu_2(B) \int_B
\vp^2(|x|)\,dx= N \mu_2(B).\nonumber
\end{align}

\begin{lem}\label{lem:appest}
In the above setting, we have
\be \label{eq:appest}
 \mu_2(U_r,g_r)\leq \frac{\displaystyle \int_{U_r}\left((G')^2+(N-1)\frac{G^2}{|x|^2}\right)dv_{g_r}
   +\displaystyle\frac{\rho_r^2}{3}\int_{U_r}\frac{G^2}{|x|^2}R^{r}_{lk}x^lx^k\,dv_{g_r}    }
   { \displaystyle  \int_{U_r}G^2\,dv_{g_r}}+O(r\rho_r^2).
\ee
as $r \to 0$. Moreover,
\begin{equation}
  \label{eq:7}
 (1+O(r^2)) \mu_2(U_r,g_r)\leq \mu_2(U_r)
\end{equation}
and
\be\label{eq:estintUrgsq}
 \int_{U_r}G^2\,dv_{g_r}\geq N-\frac{|B|\rho^2_r}{6}S(y_0)\int_0^1\vp^2t^{N+1}dt+
O(r\rho_r^2).
\ee
\end{lem}
\begin{pf}
We start by proving \eqref{eq:estintUrgsq}. Clearly
$$
\int_{U_r}G^2\,dv_{g_r}=\int_{B}G^2\,dv_{g_r}+
\int_{U_r\setminus(U_r\cap B)}G^2\,dv_{g_r}-\int_{B\setminus(U_r\cap
B)}G^2\,dv_{g_r}.
$$
Using \eqref{eq:volU}, the fact that $G$ is non-decreasing and that
$G=\vp(1)$ on $U_r\setminus(U_r\cap B)$ we get
$$
\int_{U_r}G^2\,dv_{g_r}\geq \int_{B}G^2\,dv_{g_r}+\vp(1)^2\left(
|U_r\setminus(U_r\cap B)|_{g_r}-|B\setminus(U_r\cap
B)|_{g_r}\right)=\int_{B}G^2\,dv_{g_r}.
$$
Now by \eqref{eq:metexp} and \eqref{eq:normvpxiomx} we have
\begin{align*}
 \int_{B}G^2\,dv_{g_r}&=\int_{B}\vp^2(|x|)\,dx-\frac{\rho_r^2}{6}\int_0^1\vp^2t^{N+1}dt\int_{\de B}R^{r}_{lk}x^l
 x^k d\s+O(\rho_r^3)\\
&=N-\frac{|B|\rho_r^2}{6}S(y_0)\int_0^1\vp^2t^{N+1}dt+O(r\rho_r^2).
\end{align*}
From this we conclude
$$
 \int_{U_r}G^2\,dv_{g_r}\geq N-\frac{|B|\rho^2_r}{6}S(y_0)\int_0^1\vp^2t^{N+1}dt+
O(r\rho_r^2),
$$
so that \eqref{eq:estintUrgsq} holds. Moreover, by \eqref{eq:metexp} we have
$$
|\n f_i|_{g_r}^2= |\n f_i|^2-\frac{\rho^2_r}{3}
R^{r}_{jklm}x^j\frac{\de f_i}{\de x^k} x^l\frac{\de f_i}{\de x^m}+O\left(\rho_r^3|x|^3|\n f_i|^2 \right),
$$
in $B^r$. Using furthermore that, by general properties of the
Riemannian curvature tensor, $R^{r}_{jklm}x^j x^k=0$ for every
$l,m$ and $R^{r}_{jklm} x^l x^m=0$ for every $j,k$, we obtain
$$
R^{r}_{jklm}x^j\frac{\de f_i}{\de x^k} x^l\frac{\de f_i}{\de x^m}=
\frac{G^2}{|x|^2} R^r_{jili} x^j x^l
$$
by (\ref{eq:18}). Therefore, summing over $i$ and using \eqref{smfi2},
we find that
\be\label{eq:smnfiwi} \sum_{i=1}^N|\n
f_i|_{g_r}^2=(G')^2+(N-1)\frac{G^2}{|x|^2}+
\frac{\rho_r^2G^2}{3|x|^2}R^{r}_{jl}x^jx^l+O\left(\rho_r^3|x|^3
\left( (G')^2+(N-1)\frac{G^2}{|x|^2}\right)\right).
\ee
To estimate the last term in \eqref{eq:smnfiwi},  we first note that,
since $G' \equiv 0$ in $U_r\setminus(U_r\cap B)\subset B^r$,
\begin{align*}
\rho_r^3\int_{ U_r\setminus(U_r\cap B) }|x|^3 \left( (G')^2+(N-1)\frac{G^2}{|x|^2}\right)dv_{g_r}
&=\rho_r^3G(1)^2\int_{ U_r\setminus(U_r\cap B)
}\frac{|x|^3}{|x|^2}dv_{g_r}\\
&\leq r \rho_r^2G(1)^2 |U|_{g_r}=O(r\rho_r^2).
\end{align*}
Moreover, since $|x| \le 1$ in $B$ we have
$$
\rho_r^3\int_{  B }|x|^3 \left(
  (G')^2+(N-1)\frac{G^2}{|x|^2}\right)dv_{g_r}=O(\rho_r^3)=O(r\rho_r^2).
$$
Hence we deduce that
\be\label{eq:intr32}
\int_{U_r}O\left(\rho_r^3|x|^3 \left( (G')^2+(N-1)\frac{G^2}{|x|^2}\right) \right)dv_{g_r}= O(r\rho_r^2).
\ee
Now \eqref{eq:appest} follows immediately from
\eqref{eq:mudlesSm},~\eqref{smfi2},~\eqref{eq:estintUrgsq},~\eqref{eq:smnfiwi} and
\eqref{eq:intr32}. Combining \eqref{eq:appest} and \eqref{eq:estintUrgsq}, we also deduce that
\begin{equation}
  \label{eq:6}
\mu_2(U_r,g_r) \le C+O(r^2) \qquad \text{as $r \to 0$ with a constant
  $C>0$.}
\end{equation}
Now to prove (\ref{eq:7}), we consider a normalized eigenfunction
$h_r$ corresponding to $\mu_2(U_r)$, i.e. $h_r \in H^1(U_r)$ satisfies
$$
\int_{U_r} h_r^2\,dx = 1,\quad \int_{U_r} h_r\,dx = 0 \qquad
\text{and}\qquad
\int_{U_r} |\nabla h_r|^2\,dx = \mu_2(U_r).
$$
By \eqref{eq:metexp-1} and \eqref{eq:volUdx} we then have
$$
\Bigl|\frac{1}{|U_r|_{g_r}}\int_{U_r} h_r\,dv_{g_r}\Bigr| = O(r^2)\frac{1}{|U_r|_{g_r}} \int_{U_r}|h_r|\,dx \le
O(r^2)\frac{\sqrt{|U_r|}}{|U_r|_{g_r}} = O(r^2),
$$
With $c_r=\frac{1}{|U_r|_{g_r}} \int_{U_r} h_r\,dv_{g_r}$ we therefore deduce
$$
\int_{U_r} (h_r-c_r)^2 \,dv_{g_r}= \int_{U_r} (h_r +O(r^2))^2
(1+O(r^2))\,dx= 1 +O(r^2).
$$
Therefore the variational characterization of $\mu_2(U_r,g_r)$ yields
$$
\mu_2(U_r,g_r) \le \frac{1}{1+O(r^2)} \int_{U_r} |\n h_r|_{g_r}^2
\,dv_{g_r}= (1+O(r^2)) \int_{U_r} |\n h_r|^2(1+O(r^2))\,dx
=(1+O(r^2))\mu_2(U_r),
$$
and (\ref{eq:7}) follows.
\end{pf}\\
The following lemma controls the symmetric distance
between $B$ and $U_r$ with the help of a recent stability estimate of
Brasco and Pratelli \cite{BP} for $\mu_2$ in the euclidean setting.
\begin{lem}\label{lem:asymd}
Assume that $\mu_2(U_r,g_r) \ge \mu_2(B)(1+o(1))$ as $r \to 0$ for the
family of domains $U_r$ defined in \eqref{eq:defUr}, and let $U_r
\bigtriangleup B= (U_r \cup B) \setminus (U_r \cap B)$. Then 
\be
\label{eq:estsymd}
|U_r \bigtriangleup B | \to 0 \qquad \text{as $r \to 0$.}
\ee
\end{lem}
\begin{pf}
We consider the rescaled set $U'_r=(1+\d(r))U_r$, where
$\d(r)$ is chosen such that $|U'_r|=|B|$. Then $\d(r)=O(r^2)$ by
\eqref{eq:volUdx}. By (\ref{eq:7})
and by assumption, we see that
$$
\mu_2(U_r')=(1+\d(r))^{-2}\mu_2(U_r) \ge
(1+O(r^2))\mu_2(U_r,g_r) \ge
\mu_2(B)(1+o(1)) \qquad \text{as $r \to 0$,} 
$$
whereas $\mu_2(U_r') \le \mu(B)$ by Weinberger's result \cite{Wein}.
By \cite[Theorem 4.1]{BP}, there exist points $x_r
\in \R^N$ such that 
\begin{equation}
  \label{eq:12}
 |U_r' \bigtriangleup B(x_r) |^2\leq
C(\mu_2(B)-\mu_2(U_r')) \to 0 \qquad \text{as $r \to 0$}
\end{equation}
with some constant $C>0$, where $B(x_r)$ stands for the ball in $\R^N$ centered at $x_r$ with
radius $1$. Since $\d(r)=O(r^2)$, it is easy
to see that
\begin{equation}
  \label{eq:13}
\lim_{r \to 0}|U_r \bigtriangleup
B(x_r) |= \lim_{r \to 0}|U_r' \bigtriangleup B(x_r) |= 0.
\end{equation}
Consequently, \eqref{eq:estsymd} follows once we have shown
that $x_r \to 0$ as $r \to 0$. So we
suppose by contradiction that, after passing to a subsequence,
$\inf_{r} |x_r|>0$ and $\frac{x_r}{|x_r|} \to x_0$ as $r \to 0$
for some $x_0 \in \R^N$ with $|x_0|=1$. From \eqref{eq:choOr},
\eqref{eq:metexp-1} and (\ref{eq:13}), we then infer that 
$$
\int_{B}G(|x+x_r|)\frac{x+x_r}{|x+x_r|} \cdot x_0\,dx=
\int_{B(x_r)}G(|x|) \frac{x}{|x|} \cdot x_0\,dx =
\int_{U_r}G(|x|) \frac{x}{|x|} \cdot x_0\,dx+o(1) \to 0 
$$
as $r \to 0$. If $|x_r| \to \infty$ for a subsequence, it
would follow by the definition of $G$ that 
$$
\int_{B}\frac{x+x_r}{|x+x_r|} \cdot x_0\,dx  \to 0 \qquad \text{as $r \to 0$,} 
$$
whereas, on the other hand, $\frac{x+x_r}{|x+x_r|} \to x_0$ uniformly on $B$. This is
  impossible, so we conclude that the sequence $x_r$ is bounded and
  therefore, along a subsequence, $x_r \to \tilde x \not = 0$ as $r \to
  0$ for some $\tilde x \in \R^N \setminus \{0\}$. Using \eqref{eq:choOr},
\eqref{eq:metexp-1} and (\ref{eq:13}) similarly as before, we now infer that 
\begin{equation}
  \label{eq:35}
\int_{B}G(|x+\tilde x|)\frac{x+\tilde x}{|x+\tilde x|} \cdot \tilde x\,dx = 0  
\end{equation}
Let $D:= \{x \in B\::\: x
\cdot \tilde x>0\}$, and let $\sigma: B \to B$ denote the reflection at the hyperplane $\{x \in
\R^N\::\: x \cdot \tilde x = 0\}$ given by $\sigma(x)= x-2 x \cdot
\frac{\tilde x}{|\tilde x|^2}\tilde x$. Elementary geometric considerations show that
$$
|x+\tilde x|>|\sigma(x)+\tilde x|\qquad \text{and}\qquad  
\frac{x+\tilde x}{|x+\tilde x|} \cdot \tilde x >
\Bigl| \frac{\sigma(x)+\tilde x}{|\sigma(x)+\tilde x|} \cdot \tilde x
\Bigr|\qquad \qquad \text{for $x \in D \setminus \R \tilde x$.}
$$
 Since $G(|x|)$ is nondecreasing in $|x|$ and
positive for $x \not=0$, we conclude by a change of variable that 
$$
\int_{B} G(|x+\tilde x|)\frac{x+\tilde x}{|x+\tilde x|} \cdot \tilde
x\,dx = \int_{D}\Bigl[G(|x+\tilde x|)\frac{x+\tilde x}{|x+\tilde x|} \cdot \tilde
x+ G(|\sigma(x)+\tilde x|)\frac{\sigma(x)+\tilde x}{|\sigma(x)+\tilde x|} \cdot \tilde
x\Bigr] \,dx>0,
$$ 
contradicting \eqref{eq:35}. The contradiction shows that $x_r \to 0$
as $r \to 0$, which, as remarked before, yields the claim.
\end{pf}

\begin{lem}\label{lem:upbmudUr}
Assume that
\be\label{eq:symdtozlem}
|U_r \bigtriangleup B|\to 0\quad \textrm{ as } r\to0
\ee
for the family of domains $U_r$ defined in \eqref{eq:defUr}. Then
\be\label{eq:upbmtU}
 \mu_2(U_r,g_r)\leq \left( 1-\frac{N-2-|B|\vp(1)^2
}{6N \mu_2(B)}\rho_r^2\, S(y_0)+o(\rho_r^2)
 \right)\,\mu_2(B) \qquad \text{as $r \to 0$.}
\ee
\end{lem}
\begin{pf}
 We shall estimate the terms in \eqref{eq:appest} to reach the upper bound \eqref{eq:upbmtU}.
  First note that, by \eqref{eq:symdtozlem}, (\ref{eq:11}) and \eqref{eq:normvpxiomx},
  \begin{align}
\int_{U_r}\frac{G^2}{|x|^2}R^{r}_{lk}x^l x^k dv_{g_r}
=\int_{B}\frac{G^2}{|x|^2}R^{r}_{lk}x^l x^k
dv_{g_r}+o(1)&=R_{kk} \int_{B}\vp^2(|x|)
\frac{[x^k]^2}{|x|^2}\,dx+o(1) \nonumber\\
&=S(y_0)+o(1).\label{intr31}
  \end{align}
Since, as noted in \cite[p. 636]{Wein}, the mapping $|x|\mapsto
(G')^2(|x|)+(N-1)\frac{G(|x|)^2}{|x|^2}$ is non-increasing, we have by \eqref{eq:volU}
\begin{align}
\int_{U_r}\left((G')^2+(N-1)\frac{G^2}{|x|^2}\right)dv_{g_r}
&=\int_{B}\dots dv_{g_r}+\int_{U_r\setminus(U_r\cap
B)} \dots dv_{g_r}- \int_{B\setminus(U_r\cap B)}\dots dv_{g_r} \nonumber\\
&\leq
\int_{B}\left((G')^2+(N-1)\frac{G^2}{|x|^2}\right)dv_{g_r}. \label{eq:istUrB}
\end{align}
Moreover, using \eqref{eq:metexp} and (\ref{eq:19}), we compute  
\begin{align*}
\int_{B}&\left((G')^2+(N-1)\frac{G^2}{|x|^2}\right)dv_{g_r}
=
\int_{B}\left((G')^2+(N-1)\frac{G^2}{|x|^2}\right)\Bigl(1-\frac{\rho_r^2}{6}\,
R_{lk}^{r}x^lx^k +|x|^3 {O}(\rho_r^3)\Bigr)dx\\
&=\int_{B}\left((G')^2+(N-1)\frac{G^2}{|x|^2}\right)dx- \frac{\rho_r^2}{6} \int_{\de B}
R_{lk}^{r}x^lx^k \,d\sigma\,\int_0^1\left((\vp')^2+(N-1)\frac{\vp^2}{t^2}\right)t^{N+1}dt
+{O}(\rho_r^3)\\
&= N \mu_2(B) -\frac{|B|\rho^2_r}{6}S(y_0)
 \int_0^1\left((\vp')^2+(N-1)\frac{\vp^2}{t^2}\right)t^{N+1}dt+o(\rho_r^2).
\end{align*}
Notice that, by \eqref{eq:eqvp},
\begin{align*}
\int_0^1\left(
(\vp')^2+(N-1)\frac{\vp^2}{t^2}\right)t^{N+1}dt&=\frac{1}{|B|}\int_B\vp^2(|x|)
dx-\vp(1)^2+\mu_2(B)\int_0^1\vp^2t^{N+1}dt\\
&=\frac{N}{|B|}-\vp(1)^2+\mu_2(B)\int_0^1\vp^2t^{N+1}dt\\
\end{align*}
The two equalities above and \eqref{eq:istUrB}  yield
\begin{align*}
\int_{U_r}\left((G')^2+(N-1)\frac{G^2}{|x|^2}\right)dv_{g_r}&\leq
N\mu_2(B)-\frac{\rho_r^2N}{6}S(y_0)
+\frac{|B|\vp(1)^2}{6}\rho_r^2S(y_0)\\
&-\frac{|B|\rho_r^2}{6}\mu_2(B)
S(y_0)\int_0^1\vp^2t^{N+1}dt+o(\rho_r^2).
\end{align*}
Combining this with \eqref{eq:estintUrgsq} and \eqref{intr31}, we
obtain
$$
\mu_2(U_r,g_r)\leq\mu_2(B)-\frac{N-2-|B|\vp(1)^2
}{6N}\rho_r^2 S(y_0)+ o(\rho_r^2)
$$
and the proof is complete.
\end{pf}\\

\section{Proof of the main result}
\label{sec:proof-main-result}
In this section we complete the proof of
Theorem~\ref{pro:Main-rslt}. Part (i) follows immediately from
Corollary \ref{sec:expansion-mu_2-small-2}, and the lower bound in
Part (ii) is a direct consequence of Part (i). Hence it remains to prove
the upper bound in Part (ii). For this we assume by contradiction that there exists $\e_0>0$  and sequences of numbers
$r_k> 0$ and
 $v_{r_k}\in (0\,,\,\left| B_g(y_0,r_k)\right|_g)$ such
 that $r_k \to 0$ as $k \to \infty$ and
$$
SW_{B_g(y_0,r_k)}(v_{r_k})> \left(1-\left(\g_N S(y_0) -\e_0 \right)
 \left( \frac{v_{r_k}}{|B|} \right)^{\frac{2}{N}}     \right)SW_{\R^{N}}(v_{r_k}).
$$
%
%
Then  there exist regular domains $\O_{r_k}\subset B_g(y_0,r_k)$ with
$|\O_{r_k}|_g=v_{r_k}$ and such that
\be\label{eq:mutOc}
\mu_2(\O_{r_k},g)>  \left(1-\left(\g_N S(y_0) -\e_0 \right) \left(
\frac{v_{r_k}}{|B|} \right)^{\frac{2}{N}} \right)SW_{\R^{N}}(v_{r_k}).
\ee
As in Section~\ref{s:Ubm2}, we write $r$ instead of $r_k$ in the
following. We obtain $p_r\in B_g(y_0,r)$ such  that \eqref{eq:centmass} holds and
 we define $\rho_r$, $g_r$ and $U_r$ accordingly as above.
It is easy to see from \eqref{eq:mutOc} and the scale invariance of $\O\mapsto
|\O|_g^{\frac{2}{N}}\mu_2(\O,g)$ that
\begin{equation}
  \label{eq:20}
|U_r|_{g_r}^{\frac{2}{N}}\mu_2(U_r,g_r)>\left(1-\left(\g_N S(y_0)
-\e_0
\right)\left(\frac{|B_g(p_r,\rho_r)|_g}{|B|}\right)^{\frac{2}{N}}
\right)|B|^{\frac{2}{N}}\mu_2(B).
\end{equation}
By \eqref{eq:expvolBgr} and \eqref{eq:volU}, we also find 
$$
\Bigl(\frac{|B|}{|U_r|_{g_r}}\Bigr)^{\frac{2}{N}}= 
\Bigl(\frac{\rho_r^{N}
  |B|}{|B_g(p_r,\rho_r)|_g}\Bigr)^{\frac{2}{N}}= 1 + \frac{1}{3N(N+2)}S(y_0)\rho_r^2+o(\rho_r^2)
$$
and 
$$
\Bigl(\frac{|B_g(p_r,\rho_r)|_g}{|B|}\Bigr)^{\frac{2}{N}}= \rho_r^2+o(\rho_r^2).
$$
Combining this with (\ref{eq:20}), we obtain
\begin{align}
\mu_2(U_r,g_r)&> \left(1 + \Bigl(\bigl[\frac{1}{3N(N+2)}-\g_N \bigr] S(y_0)
+\e_0\Bigr) \rho_r^2+o(\rho_r^2) \right) \mu_2(B) \nonumber\\
&=\left(1 -\Bigl(\frac{N-2-|B|\vp(1)^2
}{6N \mu_2(B)} S(y_0)-\e_0\Bigr) \rho_r^2+ o(\rho_r^2)\right) \mu_2(B)   \label{eq:mutUc}
\end{align}
and in particular $\mu_2(U_r,g_r)\ge\mu_2(B)(1+o(1))$ as $r \to 0$. From this, we can apply Lemma
\ref{lem:asymd} to get that
$$ |U_r \bigtriangleup B|\to
0\quad \textrm{ as } r\to0.
$$
 Therefore by Lemma \ref{lem:upbmudUr} we get
$$
\mu_2(U_r,g_r)\leq \left( 1-\frac{N-2-|B|\vp(1)^2
}{6N \mu_2(B)}S(y_0)\,\rho_r^2+o(\rho_r^2)
 \right)\,\mu_2(B) \qquad \text{as $r \to 0$},
$$
and this contradicts \eqref{eq:mutUc}. Hence Theorem~\ref{pro:Main-rslt}
 is proved.

\end{document}